\documentclass[letterpaper, 10 pt, conference, twocolumn]{article}
\usepackage[margin=.75in]{geometry}
\usepackage{cite,amsmath,times,amssymb,graphicx}

\title{\LARGE \bf Stability of Spatially Distributed, Intersecting Aircraft Flows Under Sequential Conflict Resolution Schemes
\thanks{This work was supported in part by NASA contract NNX08AY52A and 
        Air Force contract FA9550-08-1-0375.}%
}

\begin{document}

\author{Troy~Hand
\thanks{T. Hand is with the Daniel Guggenheim School of Aerospace
        Engineering, Georgia Institute of Technology, Atlanta, GA 30332 USA
        (e-mail: tshand@gatech.edu).} %
\and Zhi-Hong~Mao
\thanks{Z.-H. Mao is with the Department of Electrical and 
        Computer Engineering and the Department of Bioengineering, 
        University of Pittsburgh, Pittsburgh, 
        PA 15261 USA (email: maozh@engr.pitt.edu).}%
\and Eric~Feron
\thanks{E. Feron is with the Daniel Guggenheim School of Aerospace
        Engineering, Georgia Institute of Technology, Atlanta, GA 30332 USA
        (e-mail: feron@gatech.edu).}
}

\date{}
\maketitle
\thispagestyle{empty}
\pagestyle{empty}

\begin{abstract}
This paper discusses the effect of sequential conflict resolution maneuvers of an infinite aircraft flow through a finite control volume. Aircraft flow models are utilized to simulate traffic flows and determine stability. Pseudo-random flow geometry is considered to determine airspace stability in a more random airspace, where aircraft flows are spread over a given positive width. The use of this aircraft flow model generates a more realistic flow geometry. A set of upper bounds on the maximal aircraft deviation during conflict resolution is derived. Also with this flow geometry it is proven that these bounds are not symmetric, unlike the symmetric bounds derived in previous papers for simpler flow configurations. Stability is preserved under sequential conflict resolution algorithms for all flow geometries discussed in this paper.
\end{abstract}

\section{Introduction}
The current use of centralized air traffic control to ensure aircraft separation is a safe option that has been proven over the years. The process of control is typically through the use of surveillance radars, voice radio systems, limited computer support systems, and numerous complex procedures~\cite{NASA:07}. With current air traffic control techniques, increasing air traffic volume steadily increases complexity~\cite{Wol:07} and produces drawbacks such as: system bottlenecks, indirect routing, and lack of navigation freedom for airlines~\cite{MaF:01,MFB:01}, not to mention the increased workload of the ground controllers~\cite{TrM:08}. The Federal Aviation Administration (FAA) and airlines have proposed the concept of ``Free Flight''~\cite{RTCA:95} to eliminate restrictions imposed by the current system and allow for more navigation freedom as well as direct routing. The US is further developing NextGen (Next Generation Air Transportation System) to address the challenges of increasing air traffic volume as well as limitations on operational flexibility~\cite{NASA:07,HwS:08}. Europe is also further developing SESAR (Single European Sky Air-traffic-management Research program) to address growing problems with their current air traffic operations~\cite{KyM:06}. 

Automation is a key element necessary to achieve the goals set by NextGen and allow the concept of Free Flight to be more viable~\cite{GSF:09,SiM:09,ErP:02}. A decentralized air traffic control architecture could be utilized, resulting in the automation of the air traffic and possibly alleviating some of the drawbacks associated with centralized air traffic control. A decentralized solution would require each aircraft to determine its maneuver based on information shared between aircraft, such as position and velocity supplied by the Global Positioning System (GPS). A switch from current systems to GPS would alleviate limitations associated with the ground-based navigation infrastructure and lead to Free Flight. The use of a decentralized air traffic control would also distribute the work load and allow for an almost fully automated system, allowing human controllers to manage considerably more aircraft. The use of decentralized air traffic control architecture will allow the system to be scalable to increasing air traffic volume. With indications of a significant increase in air traffic volume, ranging from a factor of two to three by 2025~\cite{NASA:07,Jen:07}; decentralized air traffic control needs to be considered.

The purpose of this paper is to help build a strong analytical base for understanding conflict resolution and its limits when sequential control is utilized. Three different aircraft flow models are considered in this paper and air traffic stability is determined for each. An aircraft flow is defined as being stable if all conflicts are resolved and the conflict resolution maneuver bounded~\cite{MFB:01}. Two models are recreated from the work in ~\cite{MaF:01} to validate algorithms against previous studies. These two simpler flows are the orthogonal flow geometry and the arbitrary encounter angle geometry. Both of which are stable and have an analytical solution for the displacement bounds~\cite{MFB:01,MaF:01}. The purpose of the arbitrary encounter angle geometry is to generalize the orthogonal flow geometry for arbitrary encounter angles. The pseudo-random flow geometry is examined to generalize the orthogonal flow geometry for arbitrary flow thickness. Simulations are conducted to determine how flow thickness will effect the aircraft flow during conflict resolution. Stability is found to be achieved by the pseudo-random flow with analytical solutions derived for asymmetric displacement bounds.

\begin{figure}[tbhp]
\begin{center}   
\includegraphics[width=\linewidth]{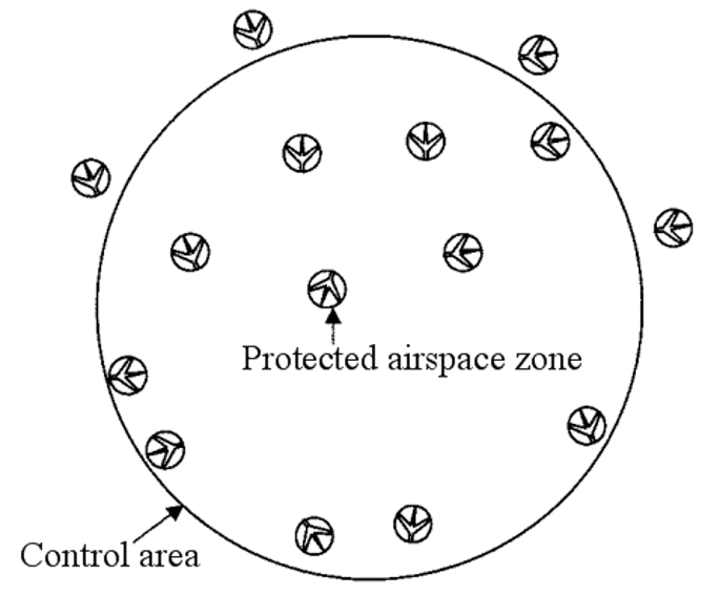} 
\caption{Aircraft flowing in and out of a control area.}\label{CA}
\end{center}
\end{figure}

\section{Problem Description}
The models presented in this paper consist of an infinite number ofl aircraft traveling sequentially through a finite airspace. Using this approach eliminates the concern about the domino effect, where one aircraft's maneuver causes another aircraft to maneuver and so on. By having an infinite sequential aircraft flow, the domino effect would be seen very clearly if it were to occur. All the models incorporate an optimized decentralized conflict resolution rule, which simply determines which conflict resolution maneuver (i.e. left or right) would result in the least displacement.

For simplifying purposes, all aircraft are assumed to fly at the same altitude and to travel at a constant velocity. With all aircraft flying at the same altitude one dimension is removed and a control area is considered for the zone of conflicts as seen in Fig.~\ref{CA}. It is assumed that the position and velocity of all aircrafts within the control area are known; this can be achieved using GPS in a context\begin{figure}[bthp]
\begin{center}   
\includegraphics[width=\linewidth]{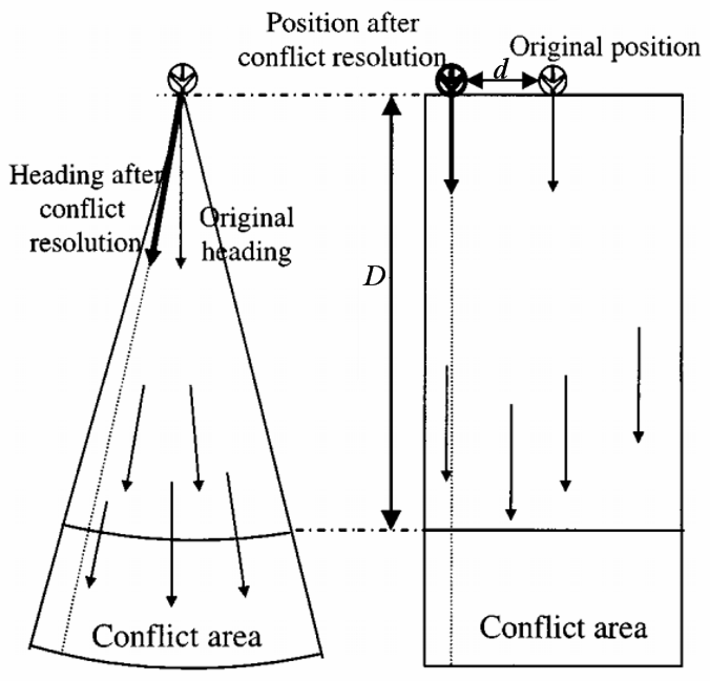} 
\caption{Heading control model $v.s.$ offset model. Left: The aircraft maneuver is an instantaneous heading change. Right: the aircraft maneuver is an instantaneous position change.}\label{man}
\end{center}
\end{figure} named ADS-B ~\cite{HwS:08}.  Therefore, the aircraft entering the control area is able to determine one maneuver that will result in the least deviation to solve conflicts with all aircraft within the control area. Once the maneuver is made, the aircraft does not make another maneuver and heading and velocity remains constant through the control area. 

An offset model is used for the conflict resolution maneuver; it is assumed to be a single lateral position change, with constant speed and heading before and after the maneuver. This model provides a close approximation to a heading change model while making analysis simpler~\cite{MaF:01,MFB:01}. While analyzing the two models in Fig.~\ref{man} it is determined that given the distance to conflict $D$, the lateral displacement $d$ in the offset model is equivalent to a heading change of amplitude $\alpha=\tan{(d/D)}^{-1}$.  If $D$ is assumed much greater than $d$ which is usually the case for strategic conflict resolution, the longitudinal displacement difference between the models is on the order of $d^2/D$, which is assumed to be small.

Three different flow geometries are considered in this paper, two of which are also considered in~\cite{MaF:01}:

\subsection{Orthogonal Flow Geometry}
\begin{figure}[tbhp]
\begin{center}   
\includegraphics[width=\linewidth]{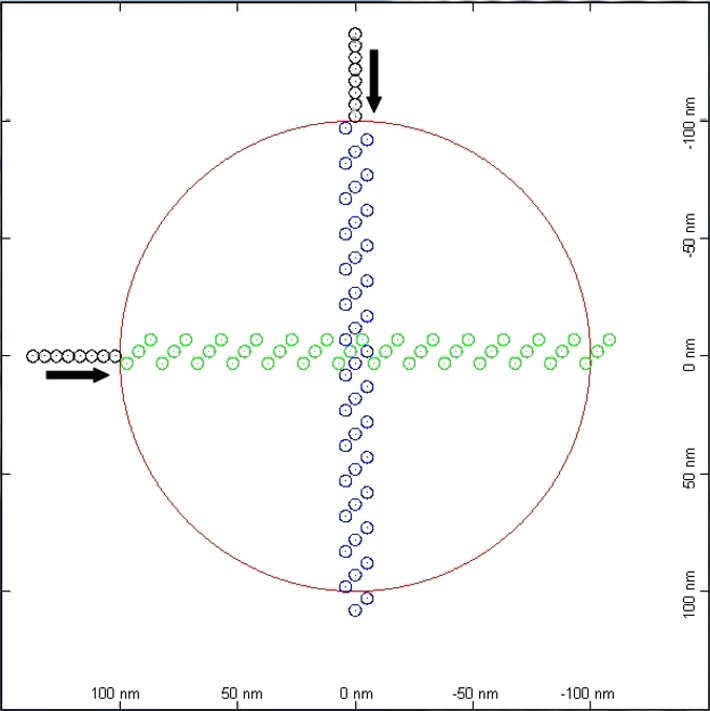} 
\caption{Orthogonal flow geometry simulation}\label{f1}
\end{center}
\end{figure}

Orthogonal flow geometry is considered for validation of optimized conflict resolutions and simulation against previous studies. This flow geometry consists of two aircraft flows; one southbound and the other eastbound. The minimum miss distance is acceptable if and only if the small circles (of radius 2.5nm) do not overlap. Aircraft in each flow maintain the same velocity; therefore, aircraft in the same flow never intersect. The aircraft in one flow needs only to consider avoiding the aircraft in the orthogonal flow by at least the miss distance specified. The control area considered for the zone of conflicts is shown in Fig.~\ref{f1} as the red circle which has a 100nm radius.

\begin{figure}[tbhp]
\begin{center}
\includegraphics[width=\linewidth]{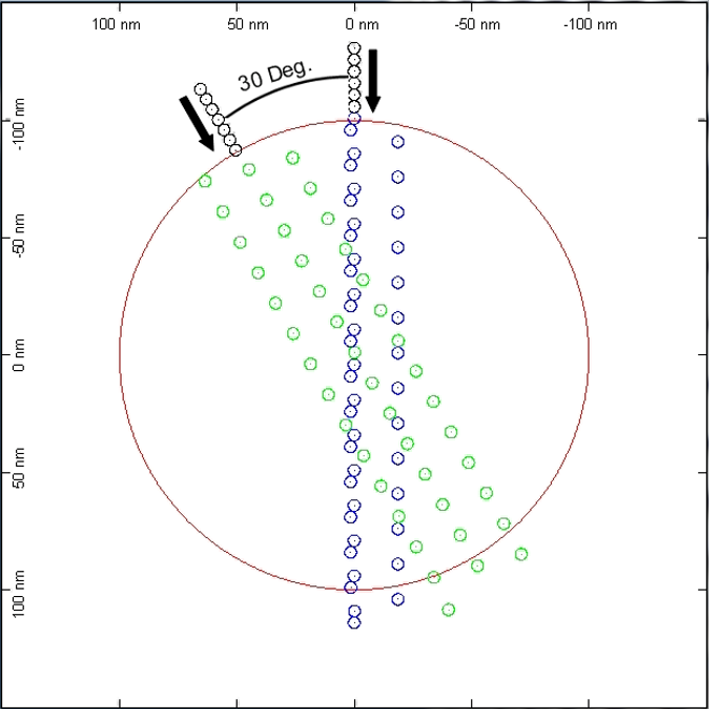}
\includegraphics[width=\linewidth]{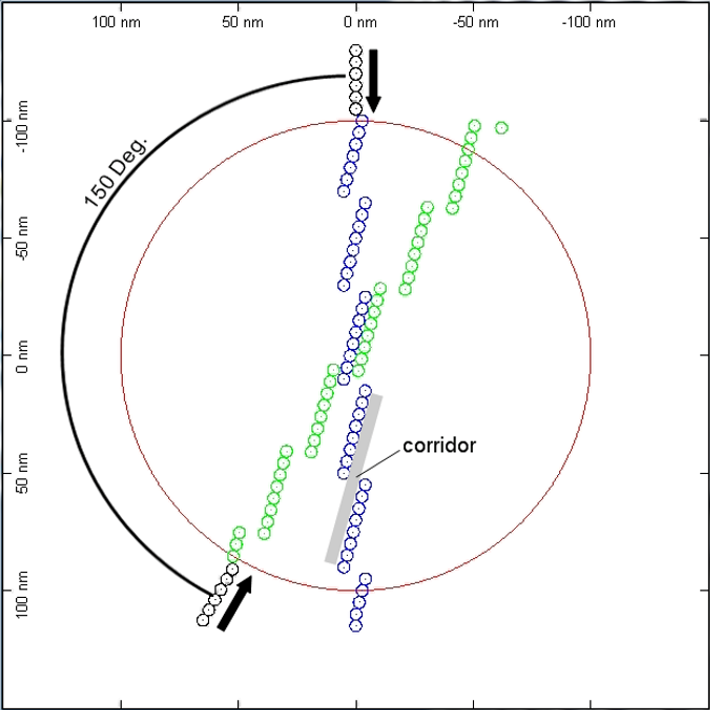}
\caption{Arbitrary encounter angle simulation}\label{f2}
\end{center}
\end{figure}

\subsection{Arbitrary Encounter Angle Geometry}
This flow geometry is a generalization of the orthogonal flow geometry for arbitrary encounter angles. This was achieved by keeping the orientation of one flow constant and angling the other flow away from the first. In the following simulations, the southbound flow has constant orientation. The eastbound flow is tested for different orientations. The encounter angle is measured counter-clockwise from the fixed flow to the other flow. The control area and miss distance are the same as those of the orthogonal flow geometry.

\subsection{Pseudo-Random Flow Geometry}
The main problem of interest for this paper is the pseudo-random flow geometry. For simplicity, the aircraft have the same speed and heading. The aircraft, however, enter the airspace at random positions within a specified starting area or ``entry gates.'' This is another generalization of the orthogonal flow geometry but for arbitrary flow thicknesses instead of arbitrary encounter angle. This flow geometry was considered in order to add disturbance to the flow yet remain close to the orthogonal flow geometry and acquire a somewhat more realistic flow.

\begin{figure}[tbhp]
\begin{center}   
\includegraphics[width=\linewidth]{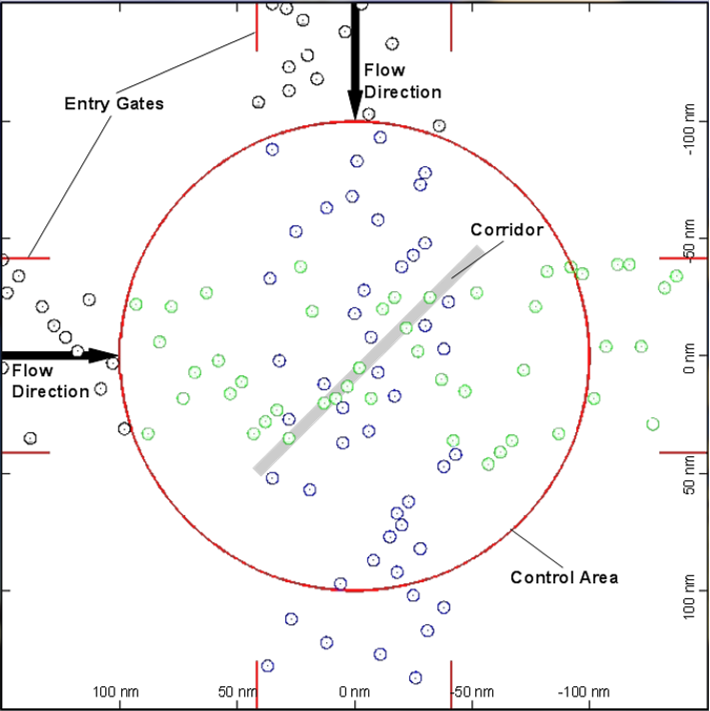} 
\caption{Pseudo-random flow geometry simulation}\label{f3}
\end{center}
\end{figure}

The control area is shown in Fig.~\ref{f3} as the red circle; and the initial flow bounds are shown as red gates located at the original flow's entrance and exit. The flow thickness is simply the distance from one entry gate to the other for a given flow. The aircraft enter sequentially but at random positions within the entry gates. Since the heading and speed is constant there are no conflicts within the same flow.

\section{Simulation}
Simulations were performed for all the flow geometries discussed in the previous section. The simulation results as seen in Fig.~\ref{f1} and Fig.~\ref{f2} match directly with the results discussed in~\cite{MaF:01,MFB:01}.  Equation 1 holds true for the bounds of the arbitrary encounter angle. Note that some aircraft take advantage of a corridor that a previous aircraft created. This corridor is the band created by propagating the safety circle of one aircraft along the constant relative velocity of the aircraft in the two flows. In Fig.~\ref{f2} notice when a small encounter angle is utilized, the offset distance the aircraft must travel for conflict resolution is greater than when the encounter angle is larger. The maximum number of aircraft involved with the same conflict increases as the encounter angle increases. It was shown in~\cite{MaF:01,MFB:01} that the maximum deviation $d_{\max}$ created by the conflict resolution maneuver satisfies:
\begin{equation}
\label{eq1}
d_{\max}=\frac{d}{\left|{\sin{(\frac{\theta}{2})}}\right|}
\end{equation}

\begin{figure}[tbhp]
\begin{center}   
\includegraphics[width=\linewidth]{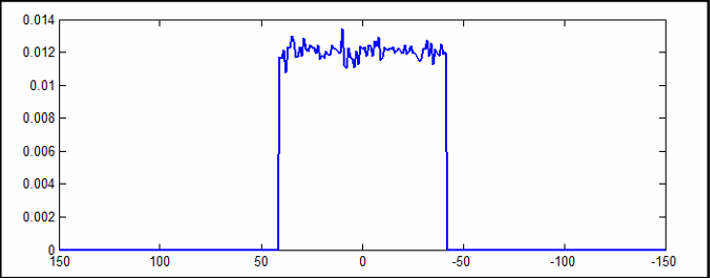} 
\caption{Exiting aircraft distribution prior to conflict resolution}\label{f4}
\end{center}
\end{figure}

\begin{figure}[tbhp]
\begin{center}   
\includegraphics[width=\linewidth]{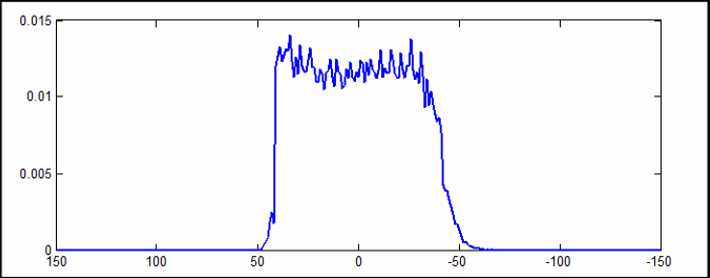} 
\caption{Exiting aircraft distribution after conflict resolution}\label{f5}
\end{center}
\end{figure}

Just as in the simpler models, the aircraft in the pseudo-random flow maneuver to fall into a corridor that a previous aircraft creates to decrease the required displacement for conflict resolution (Fig.~\ref{f3}). Also as seen in Fig.~\ref{f3} it is difficult to examine the distribution of the aircraft in any flow, and the displacement bounds are therefore harder to determine. The pseudo-random flow model was simulated for a large number of aircraft in each flow to supply sufficient data to create a probability distribution of the aircraft as they exited the airspace. This yields a clearer picture of where the aircraft will most likely end up. The graph in Fig.~\ref{f4} shows the distribution of aircraft prior to conflict resolution and Fig.~\ref{f5} shows the distribution of aircraft after conflict resolution maneuvers.

A detailed inspection of Fig.~\ref{f5} indicates that the distribution of southbound aircraft after conflict resolution has a broader support than prior to the conflict resolution. Moreover, this distribution is not symmetric, with larger possible deviations to east than to the west.


\section{Displacement Bounds }
As observed in the simulations, the displacement bounds for the pseudo-random flow geometry are asymmetric; however, we are able to show that each individual aircraft within the flow has its own displacement bounds, which are symmetric like before in the simpler model flows. Each aircraft's displacement bounds are a function of its starting position. So the bounds of the whole flow of aircrafts is simply the superposition of the displacement bounds for the individual aircrafts.

This proof aims at developing analytical solutions for the left and right displacement boundaries from the center of the aircraft flow. The first part of this proof follows fairly closely to the proof of Theorem 1 in~\cite{MFB:01}. The proof considered here focuses on bounds for the southbound aircraft flow (i.e. a displacement to the aircraft's right is west). The bounds are equally valid for the eastbound aircraft flow by symmetry of the problem.

\begin{figure}[tbhp]
\begin{center}   
\includegraphics[width=\linewidth]{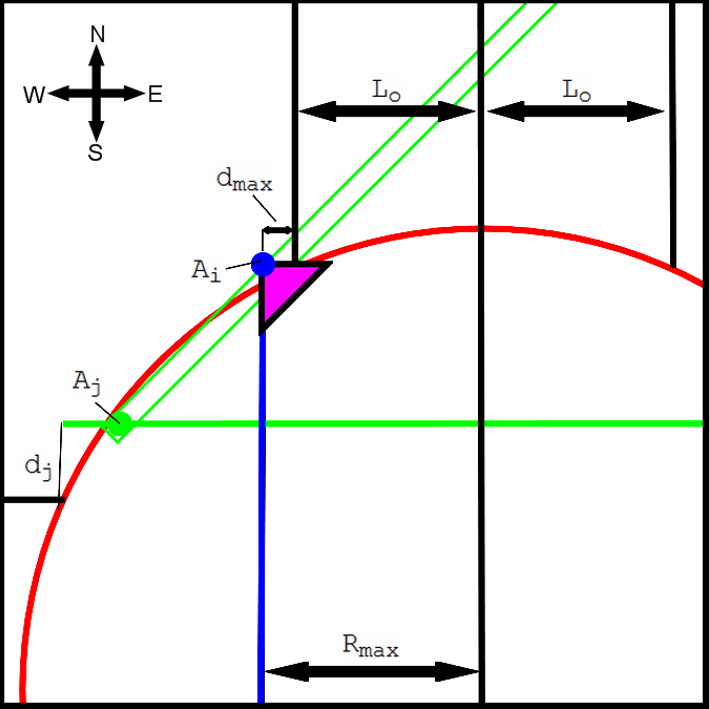} 
\caption{Pseudo-random flow geometry}\label{f6}
\end{center}
\end{figure}

First the bound for displacements to the right is found. This is accomplished by considering a southbound aircraft $A_i$ entering the control area at its far west position $L_0$, away from the center of the flow. A hypothesis is now made that states there is no maneuver of amplitude less than or equal to $d_{\max}$, where $d_{\max}=\sqrt{2}d$, in which aircraft $A_i$ is conflict free. This hypothesis implies two things: First, according to the hypothesis, the southbound aircraft $A_i$ cannot travel in a corridor created by a previous southbound aircraft whose own avoidance maneuver is small enough, because this would result in a conflict free trajectory, contradicting the hypotheses. Geometrically, no southbound aircraft can be in the pink triangular region in Fig.~\ref{f6} at the time $A_i$ makes its conflict resolution maneuver. Second, for all possible lateral deviations of $A_i$ with amplitude less than or equal to $d_{\max}$, the southbound aircraft $A_i$ is intersecting the shadow (i.e. conflict) of an eastbound aircraft, particularly $A_j$ which has already made a maneuver of amplitude $d_j$, with $d_j=2\times d_{\max}$ as shown in Fig.~\ref{f6}. A contradiction is however reached, because there are no southbound aircraft within the small pink triangular region. Therefore, the eastbound aircraft $A_j$ which was supposed to make the optimal maneuver should have had a displacement smaller than $d_j$ and both aircraft would have been conflict free. Therefore the boundary for the aircraft's displacement to the right from the center is simply the far west initial position plus the displacement $d_{\max}$. (\ref{eq2}) shows the solution to the maximum displacement to the right(west) from the center of the aircraft flow.
\begin{equation}
\label{eq2}
R_{\max}=L_0+\sqrt{2} d
\end{equation}

Because all aircraft must make an optimal maneuver, the bounds for each aircraft is symmetric, therefore the displacement an aircraft makes to the east must be strictly less than or equal to the maximum displacement to the west. The displacement bounds of all other aircrafts in the flow can be found rather easily now, since the right displacement bound is known. The right displacement bound $d_{\max}$, can be expressed in terms of $R_{\max}$ as:
\begin{equation}
\label{eq3}
d_{\max}=R_{\max}- {\textit{StartingPosition}}
\end{equation}
Where the $\textit{StartingPosition}$ is the signed distance from the center of the flow (with west being the positive direction).

By examining (\ref{eq3}) it is easy to see that an aircraft at the far east starting position would generate the largest $d_{\max}$ and therefore the largest displacement to the west. Because the aircraft must make an optimal maneuver, the largest displacement to the east is also $d_{\max}$. The bound for the aircraft's displacement to the left from the center of the flow is $2\times \textit{StartingPosition} - R_{\max}$. Thus the largest displacement to the left, computed from the center on the flow, is obtained by setting $\textit{StartingPosition}= -L_0$, thus yielding the maximum displacement to the left(east) from the center of the aircraft flow.
\begin{equation}
\label{eq4}
L_{\max}=3L_0+\sqrt{2} d
\end{equation}

Since this is a generalization of the orthogonal flow geometry the bounds can be tested for that case to determine if the same results are obtained. Since the original flow geometry as given in~\cite{MFB:01} has no flow thickness (i.e. $L_0=0$) the displacement bounds should be symmetric and equal to $d_{\max}$ obtained from (\ref{eq1}). By plugging in $L_0=0$ into (\ref{eq2}) and (\ref{eq4}) and $\theta=90^o$ into (\ref{eq1}) it is shown that $R_{\max}=L_{\max}=d_{\max}=\sqrt{2} d$.

\section{Conclusion}
With the dramatic increase in air traffic demand by 2025, a new solution to air traffic control must be considered and analytical guarantees on air traffic must be made available to guarantee system safety. The simulations demonstrated how a sequential air traffic control scheme would effect different aircraft flows. Aircraft follow the corridor created by the previous aircraft to minimize lateral displacement needed for conflict resolution. The pseudo-random flow geometry shows that as ``thickness'' is added to the flow, the lateral displacement bounds increase and become asymmetric. Although the displacement bounds increase and become asymmetric, numerical simulations show that most of the aircraft still remain within the original flow width.

\nocite{SVF:09}
\bibliography{AccSubmission_Initial}

\end{document}